\documentclass{article}
\usepackage{amssymb,amsmath}

\newcommand{\Z}{\ensuremath{\mathbf Z}}
\newcommand{\N}{\ensuremath{ \mathbf N }}
\newtheorem{theorem}{Theorem}
\newtheorem{lemma}{Lemma}
\newtheorem{corollary}{Corollary}
\newcommand{\bt}{\begin{theorem}}
\newcommand{\et}{\end{theorem}}
\newcommand{\bl}{\begin{lemma}}
\newcommand{\el}{\end{lemma}}
\newcommand{\bc}{\begin{corollary}}
\newcommand{\ec}{\end{corollary}}
\newcommand{\pf}{{\bf Proof}.\ }
\newcommand{\bq}{\begin{align*}}
\newcommand{\eq}{\end{align*}}
\newcommand{\beq}{\begin{equation}}
\newcommand{\eeq}{\end{equation}}
\newcommand{\benum}{\begin{enumerate}}
\newcommand{\eenum}{\end{enumerate}}
\newcommand{\ba}{\begin{array}}
\newcommand{\ea}{\end{array}}
\newcommand{\card}{\mbox{card}}
\newcommand{\eop}{$\square$\vspace{0.5cm}}

\begin{document}
\title{Representation functions of additive bases \\
for abelian semigroups
}
\author{Melvyn B. Nathanson\thanks{This work was supported
in part by grants from the NSA Mathematical Sciences Program
and the PSC-CUNY Research Award Program.}\\
Department of Mathematics\\
Lehman College (CUNY)\\
Bronx, New York 10468\\
Email: nathansn@alpha.lehman.cuny.edu}
\date{}

\maketitle

\begin{abstract}
Let  $X = S \oplus G,$ where $S$ is a countable 
abelian semigroup such that $S+S = S$, 
and $G$ is a countably infinite abelian group 
such that $\{ 2g : g \in G\}$ is infinite.
Let $\pi: X\rightarrow G$ be the projection map defined by 
$\pi(s,g) = g$ for all $x=(s,g) \in X.$  Let $f:X \rightarrow \N_0 \cup \{\infty\}$ 
be any map such that the set $\pi\left(f^{-1}(0)\right)$ is a finite subset 
of $G$.  Then there exists a set $B \subseteq X$ such that  
$\hat{r}_B(x) = f(x)$ for all $x \in X,$ where the restricted representation function 
$\hat{r}_B(x)$ counts the number of sets $\{x',x''\} \subseteq B$
such that $x' \neq x''$ and $x'+x''=x.$
In particular, every function $f$ 
from the integers \Z\ into $\N_0 \cup \{\infty\}$ 
such that $f^{-1}(0)$ is finite is the representation function 
of an asymptotic basis for $\Z.$
\end{abstract}

\section{Additive bases for semigroups}

Let $\N, \N_0,$ and \Z\ denote the positive integers, nonnegative integers,
and integers, respectively.
Let $X$ be an abelian semigroup, written additively,
and let $A$ and $B$ be subsets of $X$.
We define the {\em sumset}
\[
A+B = \{a + b: a\in A \mbox{ and } b \in B\}
\]
and the {\em restricted sumset}
\[
A\hat{+}B = \{a + b: a\in A, b \in B, \mbox{ and } a \neq b. \}.
\]
In particular, for $A = B$, we have 
\[
2B = B+B = \{b + b': b, b'\in B\}
\]
and
\[
2\wedge B = B \hat{+} B =  \{b + b': b, b'\in B \mbox{ and } b\neq b'\}.
\]
For every positive integer $h$ we introduce the {\em dilation}
\[
h\ast B = \{ hb : b\in B\} = \{ \underbrace{b+ \cdots+b}_{\mbox{$h$ summands}} : b\in B\}.
\]
Then
\[
2B = \left(2\wedge B\right) \cup \left( 2\ast B\right).
\]
Additive number theory is classically the study of sums of subsets 
of the semigroup of nonnegative integers.
In this paper we extend the classical theory to 
a large class of abelian semigroups. 

If $X$ is a group and $A, B \subseteq X,$ 
then we can also define the {\em difference set}
\[
A-B = \{a - b: a\in A \mbox{ and } b \in B\}.
\]
We have $-A = \{0\} - A =  \{-a:a\in A\}.$

Let $B$ be a subset of the semigroup $X$.
We associate to $B$ two representation functions,
$\hat{r}_B(x)$ and $r_B(x).$
The {\em restricted representation function}
\[
\hat{r}_B:X \rightarrow \N_0 \cup \{ \infty\}
\]
of the set $B$ counts the number of ways an element $x\in X$ can be written
as a sum of two distinct elements of $B$, that is,
\[
\hat{r}_B(x) = \card\left( \{b,b'\}\subseteq B : b+b' = x  \mbox{ and } b\neq b'\right).
\]
The {\em representation function}
\[
r_B:X \rightarrow \N_0 \cup \{ \infty\}
\]
of the set $B$ counts the number of ways an element $x\in X$ can be written
as a sum of two not necessarily distinct elements of $B$, that is,
\[
r_B(x) = \card\left( \{b,b'\}\subseteq B : b+b' = x\right).
\]

If every element of $X$ can be represented as the sum of two distinct elements of $B$, 
that is, if $2\wedge B = X$, or, equivalently, if $\hat{r}_B(x) \geq 1$ for all $x \in X$,
then the set $B$ is called a {\em restricted basis} for $X$.
If $\hat{r}_B(x) \geq 1$ for all but finitely many $x \in X$,
then $B$ is called a {\em restricted asymptotic basis} for $X$.  

Similarly, if every element of $X$ can be represented as the sum of two 
not necessarily distinct elements of $B$, that is, 
if $2B = X$, or, equivalently, if $r_B(x) \geq 1$ for all $x \in X$,
then the set $B$ is called a {\em basis} for $X$.  
If $r_B(x) \geq 1$ for all but finitely many $x \in X$,
then the set $B$ is called an {\em asymptotic basis} for $X$.  

A famous conjecture of Erd\H os and Tur\' an~\cite{erdo-tura41} 
in additive number theory states that if a set $B$ 
of nonnegative integers is an asymptotic basis or a restricted asymptotic basis 
for $\N_0$, then the representation functions
$r_B(x)$ and $\hat{r}_B(x)$ must be unbounded.
This is still an unsolved problem for the semigroup 
of nonnegative integers under addition, 
but analogues of the Erd\H os-Tur\' an 
conjecture do hold in some other abelian semigroups.  
For example, let $a$ and $b$ be nonnegative integers
and define $a \star  b = \max(a,b).$  Then $(\N_0, \star)$ is an abelian semigroup
with identity 0.
If $B$ is a nonnempty subset of $\N_0$, 
then $2B = B$ and $2\wedge B = B\setminus\{\min(B)\}$ in $(\N_0, \star)$.
It follows that if $B$ is an asymptotic basis 
or a restricted asymptotic basis for $\N_0,$
then $B$ must contain all but finitely many nonnegative integers.
Moreover, if $\card(\N_0 \setminus B) = t,$ 
then $r_B(n) = n+1-t$ and $\hat{r}_B(n) = n-t$ 
for all sufficiently large $n$, and so
\[
\limsup_{n\rightarrow \infty} r_B(n) = 
\limsup_{n\rightarrow \infty} \hat{r}_B(n) = \infty.
\]

For the semigroup $(\N,\cdot )$ of positive integers under ordinary multiplication,
Erd\H os~\cite{erdo64a} proved that if $B$ is an asymptotic basis for the multiplicative
semigroup \N, 
then the representation function $r_B(n)$ is unbounded.
Ne\u set\u ril and R\"odl~\cite{nese-rodl85} gave a simple, 
Ramsey-theoretic proof of this result,
and Nathanson~\cite{nath87} and P\u us~\cite{pus92} generalized Erd\H os' theorem 
in different directions.

The story is very different for the abelian group \Z\ of integers.
There exist subsets $B$ of integers such that $2B = \Z$ and the representation
function $r_B(x) $ is bounded.  Indeed, Nathanson~\cite{nath03}
has constructed ``dense'' bases $B$ for the integers 
such that $r_B(x) = 1$ for all $x \in \Z.$
Similarly, there exist subsets $B$ of integers such that $2\wedge B = \Z$ 
and $\hat{r}_B(x) = 1$ for all $x \in \Z$.  
A special case of a theorem in this paper is that there is no constraint
on the restricted representation functions of restricted asymptotic bases for the integers,
nor on the representation functions of asymptotic bases for the integers.
This means that for any function
\[
f:\Z \rightarrow \N_0 \cup \{\infty\}
\]
such that
\[
\card\left(f^{-1}(0)\right) < \infty,
\]
there exists a set $B \subseteq \Z$ such that $\hat{r}_B(x) = f(x)$
for all $x \in \Z.$  Indeed, we prove that if $G$ is any countable abelian group
such that $2\ast G$ is infinite, then every function
\[
f:G \rightarrow \N_0 \cup \{\infty\}
\]
with
\[
\card\left(f^{-1}(0)\right) < \infty,
\]
is the restricted representation function $\hat{r}_B$ of some set $B \subseteq G,$
and is also the representation function $\hat{r}_{B'}$ of some set $B' \subseteq G.$

In this paper we study the more general case 
of additive abelian semigroups of the form
$X = S \oplus G,$ where $S$ is a countable abelian semigroup 
and $G$ is a countably infinite abelian group.
A semigroup of the form $S \oplus G$, where $S$ is a semigroup and 
$G$ is a group, will be called a {\em semigroup with a group component.}
The main result of this paper states that 
if $X = S \oplus G$ is an abelian semigroup with a group component $G$
such  that $S+S=S$ and $2 \ast G$ is infinite, then essentially 
every function is the representation function of a
restricted asymptotic basis for $X$.
A special case of this theorem for groups was obtained by P\u us~\cite{pus91}.

\section{Representation functions of semigroups with a group component}
Let $S$ be a countable abelian semigroup, 
and let $G$ be a countably infinite abelian group.  
Both $S$ and $G$ are written additively.
We do not assume that $S$ is infinite.
We shall describe the representation functions of 
asymptotic bases and restricted 
asymptotic bases for the additive abelian semigroup
\[
X = S \oplus G = \{(s,g): s\in S \mbox{ and } g \in G\}.
\]
For $x = (s,g) \in X,$ we have the projection map
$\pi:X\rightarrow G$ defined by $\pi(x) = g$.

We begin with a simple lemma about abelian groups with finite dilation.

\bl  \label{rf:lemma:boundedorder}
Let $h \geq 2$ and let $G$ be a countably infinite abelian group.
The dilation $h\ast G$ is finite if and only if
\beq   \label{rf:group}
G \cong G_0 \oplus \left( \bigoplus_{d|h \atop d \geq 2} G_d\right) ,
\eeq
where $G_0$ is a finite abelian group 
and $G_d$ is a direct sum of cyclic groups of order $d$. 
\el

\pf
Let $G$ be a group of the form~(\ref{rf:group}).
If the positive integer $d$ divides $h$ and if $\Gamma_d$ is a cyclic group
of order $d$, then $h \ast \Gamma_d = \{0\}$ 
and so $h \ast G_d = \{0\}$.
It follows that $h\ast G = h\ast G_0$ is finite.

Conversely, suppose that $h\ast G$ is finite.
Let $\gamma \in G$.  For $k \geq 1$ we have $h^k\gamma = h(h^{k-1}\gamma),$
and so
\[
\{h^k\gamma : k = 1,2,3,\ldots\} \subseteq h\ast G.
\]
Since $h\ast G$ is finite, 
there exist positive integers $j < k$ such that
$h^j\gamma = h^k \gamma$, 
and so $\gamma$ has finite order.  Therefore, $G$ is a torsion group.

Let $m$ be the least common multiple of the orders
of the elements of the finite set $h\ast G$.
Then $(mh)\gamma = m(h\gamma) = 0$ for all $\gamma \in G,$
and so $G$ is a group of bounded order.
Since an abelian group of bounded order is a direct sum of nonnzero cyclic groups
(Kaplansky~\cite[Theorem 6]{kapl54}), we can write
\[
G = \oplus_{i=1}^{\infty} \Gamma_{d_i},
\]
where $\Gamma_d$ denotes the cyclic group of order $d$.
Let $\gamma_i$ be a generator of the group $\Gamma_{d_i}.$
Since $h\ast G$ is finite, 
it follows that $h\gamma_i = 0$ for all but finitely many $i$.
The set
\[
I_0 = \{i \geq 1: h\gamma_i \neq 0\}
\]
is finite, and
\[
G_0 = \oplus_{i\in I_0} \Gamma_{d_i}
\]
is a finite abelian group.
If $h\gamma_i = 0$, then $\gamma_i$ has order $d$ for some divisor $d$ of $h$,
$d \geq 2.$  For every divisor $d$ of $h$ with $d \geq 2,$ we define
\[
I_d = \{ i \geq 1: \mbox{$\gamma_i$ has order $d$}\}
\] 
and
\[
G_d = \oplus_{i\in I_d} \Gamma_{d}.
\]
The set $\cup_{d|h}I_d$ is infinite, and
\[
G =  G_0 \oplus \left( \oplus_{d|h}G_d \right).
\]
This completes the proof.
\eop

Let $\Gamma_2^{\infty} = \oplus_{i=1}^{\infty} \Gamma_2$
denote the direct sum of an infinite number of cyclic groups of order 2.

\bl  \label{rf:lemma:order2}
Let $G$ be a countably infinite abelian group.  Then $2\ast G$ is finite
if and only if 
\[
G \cong G_0 \oplus \Gamma_2^{\infty}
\]
for some finite abelian group $G_0$.  
\el

\pf
This is the special case of Lemma~\ref{rf:lemma:boundedorder} for $h=2.$
\eop

We consider semigroups $S$ with the property that $S+S=S.$  
Equivalently, for every $s \in S$ there exist $s', s'' \in S$ such that
$s = s' + s''.$  Every semigroup with identity has this property,
since $s = s + 0.$  If $S$ is any totally ordered set without a smallest element,
and if we define $s_1 + s_2 = \max(s_1,s_2),$ then $S$ is an abelian semigroup
such that $s = s+s$ for all $s \in S,$ but $S$ does not have an identity element.

\bt  \label{rf:theorem:1}
Let $S$ be a countable abelian semigroup such that 
for every $s \in S$ there exist $s', s'' \in S$ with $s = s' + s''.$  
Let $G$ be a countably infinite abelian group such that 
the dilation $2\ast G$ is infinite.
Consider the abelian semigroup $X = S \oplus G$  
with projection map $\pi:X\rightarrow G$.
Let 
\[
f:X \rightarrow \N_0 \cup \{\infty\}
\]
be any map such that the set 
\[
Z_0 = \pi\left(f^{-1}(0)\right)
\]
is a finite subset of $G$.  
Then there exists a set $B \subseteq X$ such that 
\[
\hat{r}_B(x) = f(x)
\]
for all $x \in X,$ where $\hat{r}_B(x)$ 
denotes the number of sets $\{b,b'\} \subseteq B$
such that $b \neq b'$ and $b+b'=x.$
\et

Note that Theorem~\ref{rf:theorem:1} is not true 
for all abelian semigroups.  For example, 
let \N\ be the additive semigroup of positive integers under addition,
and $X = \N \oplus \Z.$  For every set $B \subset X$ we have
$\hat{r}_B(1,n) = r_B(1,n) = 0$ for every $n \in \Z$.

\pf
We shall construct inductively an increasing sequence 
$B_1 \subseteq B_2 \subseteq \cdots$ of finite subsets of $X$ 
such that the set $B = \cup_{n =1}^{\infty} B_n$
has the property that $\hat{r}_B(x) = f(x)$
for all $x \in X.$

Since the dilation $2\ast G$ is infinite and the set 
$Z_0 = \pi\left(f^{-1}(0)\right)$ is finite, 
we can choose an infinite subset $U$ of $G$ such that $U \cap Z_0 = \emptyset$, 
and, if $u,u' \in U$ and $u \neq u',$ then $2u \neq 2u'.$

Let $\{x_i\}_{i=1}^{\infty}$ be a sequence of elements of $X$ such that,
for every $x \in X,$
\[
f(x) = \card\left\{i\in \N : x_i = x\right\}.
\]
Let $x_1 = (s_1,g_1) \in S\oplus G = X.$
Then $f(x_1) \geq 1.$
Since $U$ is infinite, we can choose $u_1 \in U$ such that 
\[
2u_1 \neq g_1.
\]
Choose $s'_1, s_1'' \in S$ such that $s_1 = s'_1 + s_2''.$
We define
\[
B_1 = \{(s'_1,g_1-u_1), (s''_1,u_1)\}.
\]
Then $(s_1',g_1-u_1) \neq (s_1'',u_1)$ and
\[
2\wedge B_1 = \{(s_1,g_1)\} = \{ x_1\}.
\]
Then
\[
\hat{r}_{B_1}(x) = \left\{\ba{ll} 
1 & \mbox{if $x = x_1$,}\\
0 & \mbox{otherwise,}
\ea
\right.
\]
and so
\[
\hat{r}_{B_1}(x) \leq f(x) \quad\mbox{for all $x \in X.$}
\]

Let $n \geq 2,$ and suppose that we have constructed finite sets
\[
B_1 \subseteq B_2 \subseteq \cdots \subseteq B_{n-1} \subseteq X
\]
such that 
\[
\hat{r}_{B_{n-1}}(x) \leq f(x) \quad\mbox{for all $x \in X$}
\]
and, for $i = 1,2,\ldots, n-1,$
\[
\hat{r}_{B_{n-1}}(x_i) \geq \card\left\{j \leq n-1 : x_j = x_i\right\}.
\]

Let $x_n = (s_n,g_n) \in X.$  
If $\hat{r}_{B_{n-1}}(x_n) = f(x_n),$ then we set $B_n = B_{n-1}.$ 
Suppose that $\hat{r}_{B_{n-1}}(x_n) < f(x_n).$ 
Since $U$ is an infinite subset of the group $G$
and $\pi(B_{n-1})$ and $\pi(2\wedge B_{n-1})$ are finite subsets of $G$, 
we can choose an element $u_n \in U$ satisfying the following conditions:
\benum
\item[(i)]  
\[
2u_n \neq g_n,
\]
\item[(ii)]
\[
u_n \not\in   \pi(B_{n-1})    \cup  \left( \{g_n\} - \pi(B_{n-1}) \right),
\]
\item[(iii)]
\[
u_n \not\in \pi(B_{n-1})+\{g_n\} - \pi\left(2\wedge B_{n-1} \right),
\]
\item[(iv)]
\[
u_n \not\in \pi\left(B_{n-1} \right) + \{g_n\} - Z_0,
\]
\item[(v)]
\[
u_n \not\in  \pi\left(2\wedge B_{n-1} \right) - \pi(B_{n-1}),
\]
\item[(vi)]
\[
u_n \not\in Z_0 - \pi\left( B_{n-1} \right),
\]
\item[(vii)]
\[
2u_n \not\in \pi(B_{n-1})+ \{g_n\} - \pi\left( B_{n-1} \right).
\]
\eenum
Choose $s_n', s_n'' \in S$ such that $s_n = s_n' + s_n''$, and 
let
\[
B_n = B_{n-1} \cup \{(s_n',g_n-u_n), (s_n'',u_n)\}.
\]
It follows from (i) and (ii) that
\[
(s_n',g_n-u_n) \neq (s_n'',u_n)
\]
and
\[
B_{n-1} \cap \{(s_n',g_n-u_n), (s_n'',u_n)\} = \emptyset.
\]
Therefore,
\[
2\wedge B_n = 2\wedge B_{n-1} \cup \left( B_{n-1} + \{(s_n',g_n-u_n), (s_n'',u_n)\}\right)
\cup \{x_n\}.
\]
It follows from~(iii) and~(iv) that
\[
\left( B_{n-1} + \{(s_n',g_n-u_n)\} \right) \cap 2\wedge B_{n-1} = \emptyset
\]
and that 
\[
f(x) \geq 1 \quad\mbox{for all $x \in B_{n-1} + \{(s_n',g_n-u_n)\}.$}
\]
Similarly,~(v) and~(vi) imply that
\[
\left( B_{n-1} + \{(s_n'',u_n)\} \right) \cap 2\wedge B_{n-1} = \emptyset
\]
and that 
\[
f(x) \geq 1 \quad\mbox{for all $x \in B_{n-1} + \{(s_n'',u_n)\}.$}
\]
It follows from~(vii) that
\[
\left( B_{n-1} + \{(s_n',g_n-u_n)\} \right) \cap 
\left( B_{n-1} + \{(s_n'',u_n)\} \right) = \emptyset.
\]
Conditions~(ii) and~(iii) imply that
\[
x_n \not\in B_{n-1} + \{(s_n',g_n-u_n), (s_n'',u_n)\}.
\]
Therefore,
\[
\hat{r}_{B_n}(x_n) = \hat{r}_{B_{n-1}}(x_n) + 1 
\geq \card\left\{j \leq n : x_j = x_n\right\}
\]
and
\[
\hat{r}_{B_n}(x) =
\left\{\ba{ll}
\hat{r}_{B_{n-1}}(x) & \mbox{for $x \in \left( 2\wedge B_{n-1}\right)\setminus\{x_n\}$,}\\
1 & \mbox{for $x \in B_{n-1} + \{(s_n',g_n-u_n), (s_n'',u_n)\},$}
\ea\right.
\]
hence $\hat{r}_{B_n}(x) \leq f(x)$ for all $x \in X.$

This construction produces a sequence $\{B_n\}_{n=1}^{\infty}$ of finite sets 
with the property that the infinite set $B = \cup_{n=1}^{\infty} B_n$ 
satisfies $\hat{r}_{B}(x) = f(x)$ for all $x\in X.$
This completes the proof.
\eop

\bt  \label{rf:theorem:2}
Let $G$ be a countably infinite abelian group such that 
the dilation $2\ast G$ is infinite.
Let 
\[
f:G \rightarrow \N_0 \cup \{\infty\}
\]
be any map such that the set 
\[
Z_0 = f^{-1}(0)
\]
is a finite subset of $G$.  
Then there exists a restricted asymptotic basis $B$ of order 2 for $G$ such that 
\[
\hat{r}_B(x) = f(x)
\]
for all $x \in X,$ where $\hat{r}_B(x)$ 
denotes the number of sets $\{b,b'\} \subseteq B$
such that $b \neq b'$ and $b+b'=x.$
\et

\pf
This follows immediately from Theorem~\ref{rf:theorem:1}
with $S = \{0\}.$
\eop

By the same methods we can prove the following theorems about unrestricted representation
functions.  In the unrestricted case, there are additional conditions that must be
satisfied for the proofs to work.  It suffices to replace $2\ast G$ by $12 \ast G$ 
in the statements of the Theorems~\ref{rf:theorem:3} and~\ref{rf:theorem:4}.

\bt  \label{rf:theorem:3}
Let $S$ be a countable abelian semigroup such that 
for every $s \in S$ there exist $s', s'' \in S$ with $s = s' + s''.$  
Let $G$ be a countably infinite abelian group such that 
the dilation $12\ast G$ is infinite.
Consider the abelian semigroup $X = S \oplus G$  
with projection map $\pi:X\rightarrow G$.
Let 
\[
f:X \rightarrow \N_0 \cup \{\infty\}
\]
be any map such that the set 
\[
Z_0 = f^{-1}(0)
\]
is finite. 
Then there exists a set $B \subseteq X$ such that 
\[
r_B(x) = f(x)
\]
for all $x \in X,$ where $r_B(x)$ 
denotes the number of sets $\{b,b'\} \subseteq B$
such that $b+b'=x.$
\et

\bt  \label{rf:theorem:4}
Let $G$ be a countably infinite abelian group such that 
the dilation $12\ast G$ is infinite.
Let 
\[
f:G \rightarrow \N_0 \cup \{\infty\}
\]
be any map such that the set 
\[
Z_0 = f^{-1}(0)
\]
is finite.  
Then there exists an asymptotic basis $B$ of order 2 for $G$ such that 
\[
r_B(x) = f(x)
\]
for all $x \in X,$ where $r_B(x)$ 
denotes the number of sets $\{b,b'\} \subseteq B$
such that $b+b'=x.$
\et

\section{Bases for groups of exponent 2}
Let $X = S \oplus G$ be an abelian semigroup with group component $G$.
If $2\ast G$ is finite, then it is an open problem to classify 
the representation functions of asymptotic bases 
and restricted asymptotic bases for $X$.  We know from
Lemma~\ref{rf:lemma:order2} that $G \cong G_0 \oplus \Gamma_2^{\infty},$
where $G_0$ is a finite abelian group.  Replacing the semigroup $S$
with $S \oplus G_0,$ we see that it suffices to consider semigroups of the form
$X = S \oplus \Gamma_2^{\infty}$.  Even the special case $S = \{0\}$
and $X = \Gamma_2^{\infty}$ is a mystery.  If $x,x' \in \Gamma_2^{\infty}$,
then $x + x' = 0$ if and only if $x = x',$ and so, for any subset $B$ of
$\Gamma_2^{\infty}$, we have $\hat{r}_B(0) = 0$, $r_B(0) = \card(B),$
and $\hat{r}_B(x) = r_B(x)$ for all $x \neq 0.$

The following result shows another constraint on the 
representation functions of asymptotic bases for $\Gamma_2^{\infty}$.

\bl  \label{rf:lemma:Gamma2}
Let $G = \Gamma_2^{\infty}$ be an infinite direct sum of cyclic
groups of order 2.  Let $B$ be a subset of $G$, and let $\hat{r}_B(x)$
be the restricted representation function of $B$.
If $\hat{r}_B(x) \geq 2$ for some $x \in G$, then there exist elements
$y,z \in G$ such that $\hat{r}_B(y) \geq 2$ and $\hat{r}_B(z) \geq 2$, and 
the elements $x, y,$ and $z$ are distinct.
\el

\pf
If $\hat{r}_B(x) \geq 2,$ then there exist distinct elements $a,b,c,d \in G$ such that
\[
x = a+b = c+d.
\]
Since every element of $G$ has order 2, we have
\[
y = a+c = a-c = d-b = b+d
\]
and
\[
z = a+d = a-d = c-b = b+c,
\]
and so  $\hat{r}_B(y) \geq 2$ and  $\hat{r}_B(z) \geq 2.$ 
The elements $x,y,$ and $z$ are distinct.  This completes the proof.
\eop

\bt
Let $G = \Gamma_2^{\infty}$, and let $f: G\rightarrow \N_0\cup\infty$
be a function such that $f(x) \geq 1$ for all $x \in G,$ and 
$f(x) \geq 2$ for exactly one or two elements of $G$.  
There does not exist a set
$B$ in $G$ such that $\hat{r}_B(x) = f(x).$
\et

\pf
This follows immediately from Lemma~\ref{rf:lemma:Gamma2}.
\eop

It is not hard to construct a basis $B$ for $\Gamma_2^{\infty}$
such that $r_B(x) = 1$ for all $x \neq 0$, 
but it is an open problem to describe all representation functions 
for this group.

\end{document}